\documentclass[conference]{IEEEtran}
\IEEEoverridecommandlockouts
\PassOptionsToPackage{bookmarks=false}{hyperref}
\usepackage{cite}
\usepackage{hyperref}
\usepackage{amsmath,amssymb,amsfonts, amsthm, mathtools}
\usepackage{algorithmic}
\usepackage{graphicx, subfig}
\usepackage{textcomp}
\usepackage{xcolor}
\def\BibTeX{{\rm B\kern-.05em{\sc i\kern-.025em b}\kern-.08em
    T\kern-.1667em\lower.7ex\hbox{E}\kern-.125emX}}

\newtheorem{defn}{Definition}
\newtheorem{thm}{Theorem}

\newcommand{\R}{\mathbb{R}}
\newcommand{\N}{\mathbb{N}}
\newcommand{\Part}{\mathcal {P}}
\newcommand{\Clus}{\mathcal {C}}
\newcommand{\Cover}{\mathcal {U}}
\newcommand{\Tower}{\mathfrak{U}}
\newcommand{\MM}{\mathrm{MM}}
\newcommand{\Multi}{\mathcal{M}}
\newcommand{\eps}{\varepsilon}
\newcommand{\abs}[1]{\left\lvert #1 \right\rvert}
\newcommand\restrict[1]{\raisebox{-.5ex}{$|$}_{#1}}
\makeatletter
\def\blfootnote{\xdef\@thefnmark{}\@footnotetext}
\makeatother

\begin{document}

\title{Multimapper: Data Density Sensitive Topological Visualization\\
}

\author{
\IEEEauthorblockN{Bishal Deb*\thanks{*This work was completed during internship at Adobe Systems, Noida.}}
\IEEEauthorblockA{
\textit{Chennai Mathematical Institute}\\
Siruseri, India \\
bishal@cmi.ac.in}
\and
\IEEEauthorblockN{Ankita Sarkar*}
\IEEEauthorblockA{
\textit{Chennai Mathematical Institute}\\
Siruseri, India \\
ankita\textunderscore s@cmi.ac.in}
\and
\IEEEauthorblockN{Nupur Kumari}
\IEEEauthorblockA{
\textit{Adobe Systems Inc}\\
Noida, India\\
nupurkmr9@gmail.com}
\and
\IEEEauthorblockN{Akash Rupela}
\IEEEauthorblockA{
\textit{Adobe Systems Inc}\\
Noida, India\\
arupela@adobe.com}
\and
\IEEEauthorblockN{Piyush Gupta}
\IEEEauthorblockA{
\textit{Adobe Systems Inc}\\
Noida, India\\
piygupta@adobe.com}
\and
\IEEEauthorblockN{Balaji Krishnamurthy}
\IEEEauthorblockA{
\textit{Adobe Systems Inc}\\
Noida, India\\
kbalaji@adobe.com}
}


\maketitle

\begin{abstract}\blfootnote{This paper was accepted at ICDMW}
Mapper is an algorithm that summarizes the topological information contained in a dataset and provides an insightful visualization. It takes as input a point cloud which is possibly high-dimensional, a filter function on it and an open cover on the range of the function. It returns the nerve simplicial complex of the pullback of the cover. Mapper can be considered a discrete approximation of the topological construct called Reeb space, as analysed in the $1$-dimensional case by \cite{carri19statistical}. Despite its success in obtaining insights in various fields such as in \cite{kamruzzaman2016characterizing}, Mapper is an \textit{ad hoc} technique requiring lots of parameter tuning. There is also no measure to quantify goodness of the resulting visualization, which often deviates from the Reeb space in practice. In this paper, we introduce a new cover selection scheme for data that reduces the obscuration of topological information at both the computation and visualisation steps. To achieve this, we replace global scale selection of cover with a scale selection scheme sensitive to local density of data points. We also propose a method to detect some deviations in Mapper from Reeb space via computation of persistence features on the Mapper graph.
\end{abstract}


\section{Introduction}\label{sec:intro}
Real-world data is often very high dimensional and hence not visualizable. The Mapper algorithm, developed in \cite{singh2007topological}, is a method to visualize such data in low dimensions while trying to remain true to the topological structure of data in the higher dimension. The algorithm returns a simplicial complex\footnote{\cite{Hatcher} provides a rigorous mathematical treatment of simplicial complexes.} which is a representation of data via a far less number of nodes than the number of data points. Through this visualization, it becomes convenient to gain insights from data.

%
%

The Mapper algorithm begins by applying a function $f: X \rightarrow Z$ on the input space $X$. The resulting lower dimensional image space $Z$ is covered by a set of bins that overlap each other, with every part of $Z$ included in at least one bin. Then, a clustering algorithm is applied within the $f^{-1}$ of each bin. The \emph{nerve} of these set of clusters, computed as in \ref{sec:background}, is called \emph{Mapper}. The Mapper restricted to only its nodes and edges is called the \emph{Mapper graph}.\\
Although the Mapper algorithm has been highly successful, it is difficult to work with due to a large number of parameters involved in the choices of lens function, type of cover and clustering algorithm. \cite{carri19statistical} have studied parameter selection in the case where $Z = \R$. They prove that in the $1$-dimensional setup, the Mapper graph statistically converges to a geometric structure called the Reeb graph, which encodes the topological information of the original space. This, in turn, gives a method to tune its parameters to best approximate the Reeb graph.\\Even with best parameters Mapper provides a visualization of data at a fixed scale at which the cover was constructed. Since our input space is a high dimensional discrete point cloud and not a continuous space, the best parameters still provide an approximation to its topology and it is possible that more insights are available when data is viewed at different scales. \cite{dey2016multiscale} address this issue by proposing Multiscale Mapper, where data is seen along a tower of covers. But still, each cover in the tower is at a scale which views all parts of the data through the same lens irrespective of the data distribution being different in different regions. In this paper, we propose an algorithm which allows us to view denser and sparser parts of the data at separate scales in the same visualization. This, unlike the current Mapper algorithm, is not restricted to one global scale; hence, it prevents the shattering of sparser data subsets under finer scales that are suitable for denser data subsets. Further taking inspiration from \cite{CatReeb,carri19statistical} we can characterize the output of Mapper (with lens function $f$) as \textit{good} when it approximates the Reeb space (a generalisation of the Reeb graph in higher dimensions) under $f$. 
\subsection{Contributions}
The main contributions of this paper are as follows: 
\begin{enumerate}

\item We propose \emph{Multimapper} in \ref{sec:multimapper}, which combines locally optimal Mappers into a single simplicial complex. A crucial issue with parameter selection is that the same choice of parameters might not be locally optimal for each part of data. Denser parts of the data might only reveal their detailed geometry at a finer scale, at which the sparser parts might shatter. Multimapper thus gives a more accurate representation of the data compared to any single global choice of Mapper parameters. 
\item In \ref{sec:det bad} we present a data-agnostic method of partially characterising parts of Mapper which are provably different from the corresponding Reeb graph. It will also help in identifying locally optimal scales for the cover used in Mapper. 
%
%
%
%
%
%
%
%
%

\item Lastly in \ref{sec:brickcover} we propose brick cover, a covering scheme that is computable as efficiently as the box-like cuboidal cover and produces a lower dimensional simplicial complex, i.e. at most $2$-simplices under a $2$-dimensional lens function. Thus it gives a visualisation without loss of topological information contained in the Mapper output. Under standard covering schemes, which are made of boxes i.e. $n$-fold direct product of intervals in $\R^n$, the Mapper might include higher order simplices that cannot be visualized.\\
%
%
%
%
%
%
%
%
%
%
%

\end{enumerate}

In \ref{sec:background}, we first briefly introduce the mathematical notions behind Mapper and then discuss our contributions towards improvement on Mapper algorithm. We have implemented all our experiments using $2$-dimensional lens function; however, a vast majority of it is generalizable to higher dimensions.  Since $\R$-valued lens functions limit the topological information to only edges; it is more informative to use higher dimensional lens function. But in dimensions greater than $3$, we get higher-order simplices that are hard to visualize, so we restrict ourselves to $2$ dimensions. \\
%
%


\section{Theoretical Background and Existing Work}\label{sec:background}
Below, we establish the required topological terminology.
\begin{defn}\label{def:OpenCover} An open cover $\Cover=\{U_\alpha : \alpha \in A\}$ of a space $X$ is a collection of open sets such that each point in the space is in at least one of these open sets.  
\end{defn}
In this paper, we shall refer to the individual elements of open cover as \textit{bins}. We can conceptualize covering as putting each element in one or more of these bins.

\begin{defn}\label{def:simplex}
A \emph{$k$-simplex} is the smallest convex set containing a given set of $k+1$ affinely independent points, where $u_0,u_1,\ldots,u_k$ are called affinely independent if $u_1-u_0, u_2-u_0,\ldots,u_k-u_0$ are linearly independent.
%
%
\end{defn}
\begin{defn}\label{def:face}
An $m$-simplex $\sigma'$ is said to be an \emph{$m$-face} of a $k$-simplex $\sigma$ if $m < k$ and the vertices of $\sigma'$ are a proper subset of the vertices of $\sigma$.
\end{defn}
\begin{defn}\label{def:simplicial complex}
A \emph{simplicial complex} $K$ is a set of simplices such that:
\begin{itemize}
\item A face of a simplex from $K$ is also in $K$
\item $\forall \sigma_1, \sigma_2 \in K, \sigma_1\cap\sigma_2$ is a face of both $\sigma_1$ and $\sigma_2$
\end{itemize}
\end{defn}
For an example, $1$-simplex is a line segment, $2$-simplex a triangle, $3$ simplex a tetrahedron and so on. And for a $3$-simplex (tetrahedron), $0$-faces are its vertices, $1$-faces are its edges, and $2$-faces are its triangular sides.
\begin{defn}\label{def:nerve}
Given a cover $\Cover$ of a space $X$, the nerve $N(\Cover)$ is a simplicial complex constructed as follows:
\begin{itemize}
\item The vertices (\textit{nodes}) of $N(\Cover)$ correspond to bins of $\Cover$
\item For each $k+1$ bins of $\Cover$ that have mutual non-empty intersection in $X$,  $N(\Cover)$ contains a $k$-simplex with the corresponding nodes as its vertices.
\end{itemize}
\end{defn}

The Mapper algorithm is motivated by the Nerve Theorem \cite[Corollary 4G.3]{Hatcher}, originally proposed by Pavel Alexandrov.

\begin{thm}[Nerve Theorem]\label{thm:nerve} If $\Cover$ is an open cover of a paracompact space $X$ such that every non-empty intersection of finitely many sets in $\Cover$ is contractible, then $X$ is homotopy equivalent to the nerve $N(\Cover)$.
\end{thm}
An intersection being \textit{contractible} intuitively means that we should be able to continuously shrink it to a point in $X$. Finally, \textit{homotopy equivalent} is a mathematical notion of the shape being similar -- which tells us that if the required conditions are satisfied, then the nerve would give us the shape of $X$ itself. \\

In the context of data analysis, we work with a point cloud lying in $\mathbb{R}^n\;n \in N$, which is a metric space. Every metric space is paracompact as proved in \cite{steen1978counterexamples}, hence the assumption of paracompact space is satisfied in our case. We can compute topological properties of point cloud by assuming that it is sampled from a paracompact space in $\mathbb{R}^n $ which we refer to as $X$.\\

\textbf{Mapper Algorithm}: based on these ideas, the algorithm works as follows:
\begin{enumerate}
\item Given a point cloud $X$, we project it onto a lower dimension space $Z$ by a \emph{lens function} $f$. We create a cover $\Cover$ on the image. The pre-image of each bin under $f$ then gives us a cover  $f^{-1}(\Cover)$ of $X$: 
\begin{equation}
    \begin{aligned}
        f^{-1}(\Cover) = \{f^{-1}(U_{\alpha}) : U_{\alpha} \in \Cover\}
    \end{aligned}\label{eq1}
\end{equation}

\item We obtain a modified \emph{pullback} cover $f^*(\Cover)$ of $X$ from the bins of $f^{-1}(\Cover)$. The \emph{pullback of $\Cover$ under $f$} is defined as 
\begin{equation}
f^*(\Cover) = \{C : \exists V \in f^{-1}(\Cover), C \in \Part(V)\}\label{eq2}
\end{equation}
where $\Part(V)$ is the set of path connected components of $V$. In the discrete setting, path connected components are approximated by clusters. Hence in practice, we obtain $f^*(\Cover)$ by clustering within each bin of $f^{-1}(\Cover)$.\\
\item We compute the nerve of $f^*(\Cover)$ in $X$. This nerve is the \emph{Mapper} $M(X,\Cover,f)$ of $X$.
\end{enumerate}
\vspace{3 mm}
Thus
\begin{defn}\label{def:Mapper}
Given a space $X$, a lens function $f$, and a cover $\Cover$ of $f(X)$, the \emph{Mapper} is defined as 
\begin{equation}
    M(X,\Cover,f) = N(f^*(\Cover))\label{eq3}
\end{equation}
\end{defn}
The \emph{$1$-skeleton} of Mapper is a graph consisting of simplices till 1 dimension from Mapper. Hereafter we shall omit $X, \Cover, f$ wherever they are clear from context.\\


Mapper has been applied in a wide range of usecases to get useful insights. The data analytics company Ayasdi have extensively used TDA to provide solutions in various fields \cite{AyasdiMain}, like finance \cite{AyasdiFinance}, healthcare \cite{AyasdiHealth}, sports \cite{AyasdiNBA} and machine learning \cite{AyasdiML}. \cite{kamruzzaman2016characterizing} have used Mapper to study environmental stressors on plant phenotypes. \cite{vejdemo2012topology} have used Mapper to identify voting patterns in the US House of Representatives, and have comparatively presented a marked improvement on the classification possible compared to Principal Component Analysis \cite{pearson1901principal}.\\
However, the parameters in Mapper, i.e. $f$ (lens function), $r$ (bin diameter) and $g$ (bin overlap) require tuning by trial and error, and hence involve some implicit knowledge of the data. Moreover, once the Mapper graph is obtained, insight mining is primarily done by human intervention. These are significant roadblocks in the goal of using Mapper to gain insights from truly unsupervised big data. We have made progress in automating these processes.\\
Our work results from \cite{dey2016multiscale}, who introduced and significantly developed the notion of Multiscale Mapper. This technique is based on building a \textit{tower} of covers of various scales, and studying the variations in the resulting \textit{tower} of nerves. Below, we define Multiscale Mapper and related terminology as presented in \cite{dey2016multiscale}.
\begin{defn}\label{def:tower}
A \emph{tower} of $\mathfrak{W}$ of objects with \emph{resolution} $\mathrm{res}(\mathfrak{W}) = r \in \R$ is a collection $\mathfrak{W} = \{\mathcal{W}_\eps\}_\eps\geq r$ of objects $\mathcal{W}_\eps$ together with maps $w_{\eps,\eps'} : \mathcal{W}_\eps \rightarrow \mathcal{W}'_\eps$ so that $w_{\eps,\eps'} = \mathrm{id}$ and $w_{\eps,\eps''} = w_{\eps',\eps''}\circ w_{\eps,\eps'}$ for all $r \leq \eps \leq \eps' \leq \eps''$
\end{defn}
\begin{defn}\label{def:map of covers}
Given two covers $\Cover = \{U_a\}_{a \in A}$, $\mathcal{V} = \{V_b\}_{b \in B}$, a \emph{map of covers} from $\Cover$ to $\mathcal{V}$ is a set map $\xi : A \rightarrow B$ such that $\forall a \in A, \Cover_a \subseteq \mathcal{V}_{\xi(a)}$. By abuse of notation, $\xi$ also represents the induced map $\Cover \rightarrow \mathcal{V}$.
\end{defn}
\begin{defn}\label{def:tower of covers}
A tower $\mathfrak{W}$ where the objects $\mathcal{W}_\eps$ are covers and the maps $w_{\eps,\eps'}$ are maps of covers, $\mathfrak{W}$ is called a \emph{tower of covers}.
\end{defn}
\begin{defn}\label{def:simplicial map}
Given two simplicial complexes $K,L$, a \emph{simplicial map} from $K$ to $L$ is a map $\xi : K \rightarrow L$ such that:
\begin{itemize}
\item For each vertex $v_K \in K$, $\xi(v_K) \in L$ is also a vertex
\item For each simplex $\sigma$ with vertices $\{v_0,\ldots,v_k\}$, $\xi(\sigma) \in L$ with vertices $\{\xi(v_0),\ldots,\xi(v_k)\}$.
\end{itemize}
\end{defn}
Note that if $\xi$ is not injective on vertices, then a $k$-simplex might map to a $k'$-simplex, $k' < k$.
\begin{defn}\label{def:tower of simplices}
A tower $\mathfrak{W}$ where the objects $\mathcal{W}_\eps$ are simplical complexes and the maps $w_{\eps,\eps'}$ are simplicial maps, $\mathfrak{W}$ is called a \emph{tower of simplicial complexes}.
\end{defn}
\begin{defn}\label{def:multiscale mapper}
Given a space $X$, lens function $f$, and a tower of covers $\Tower = \{\Cover_\eps\}_\eps of f(X)$, we define a \emph{Multiscale Mapper} \begin{equation}
    \MM(X,\Tower,f) = \{M(X,\Cover, f) : \Cover \in \Tower\}\label{eq4}
\end{equation}
\end{defn}
As in the case of Mapper, we omit $X, \Tower, f$ from notation wherever it is clear from context.\\
The following facts establish that the successive relationship between covers of $f(X)$ in $\Tower$ naturally corresponds to a successive relationship of Mappers within $\MM(\Tower)$.
\begin{itemize}
\item A map of covers $\xi : \Cover \rightarrow \mathcal{V}$ induces a simplicial map $N(\xi) : N(\Cover) \rightarrow N(\mathcal{V})$ by the following rule : if a vertex $u \in N(\Cover)$ corresponds to $U \in \Cover$ and a vertex $v \in N(\mathcal{V})$ corresponds to $V \in \mathcal{V}$ such that $\xi(U) = V$, then $N(\xi)(u) = v$.\\ Moreover, if $\Cover \xrightarrow{\xi} \mathcal{V} \xrightarrow{\zeta} \mathcal{W}$ are maps between covers, then $N(\zeta\circ\xi) = N(\zeta)\circ N(\xi)$. Thus a tower of covers induces a corresponding tower of simplicial complexes i.e. the nerves of each cover.
\item A map of covers $\xi : \Cover \rightarrow \mathcal{V}$ induces a map of covers between their respective pullbacks under a function $f$ \begin{equation}
f^*(\xi) : f^*(\Cover) \rightarrow f^*(\mathcal{V})\label{eq5}
\end{equation}
\begin{defn}\label{def:pullback of tower}
The \emph{pullback under $f : X \rightarrow Z$ of a tower of covers $\Tower=\{\Cover_\eps\}_\eps$ of $Z$} is defined as 
\begin{equation}f^*(\Tower) = f^*(\{\Cover_\eps\}_\eps) =\{f^*(\Cover_\eps)\}_\eps\label{eq6}
\end{equation}
\end{defn}
and is itself a valid tower of covers via the induced maps.
\end{itemize}


\section{A Selective Magnification Scheme -- Multimapper}\label{sec:multimapper}
\begin{figure*}
\centering
\subfloat[Bin size = $1/5$\textsuperscript{th} of image diameter]{\includegraphics[width=0.265\linewidth]{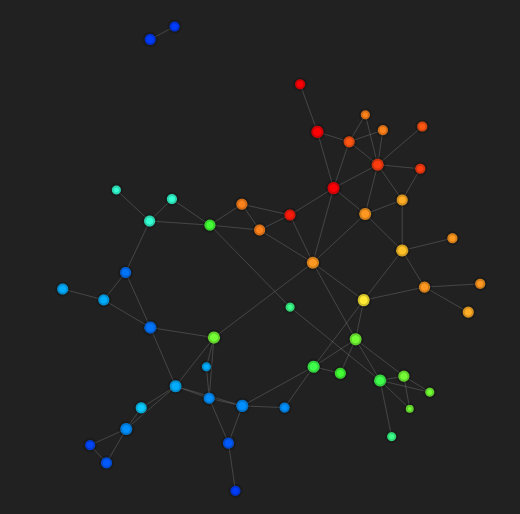}}
\qquad
\subfloat[Bin size = $1/8$\textsuperscript{th} of image diameter]{\includegraphics[width=0.27\linewidth]{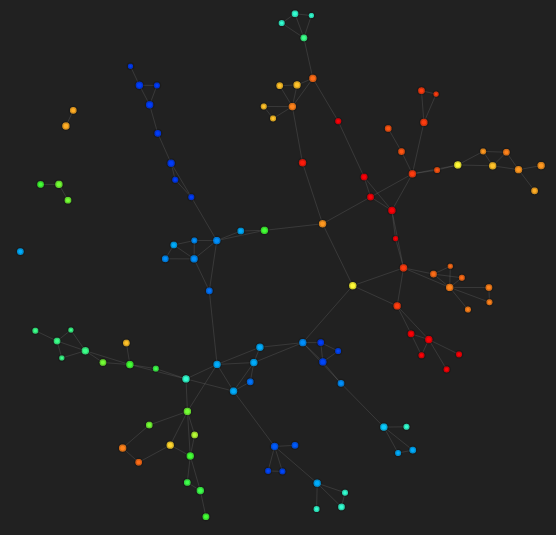}}
\qquad
\subfloat[Bin size = $1/10$\textsuperscript{th} of image diameter]{\includegraphics[width=0.27\linewidth]{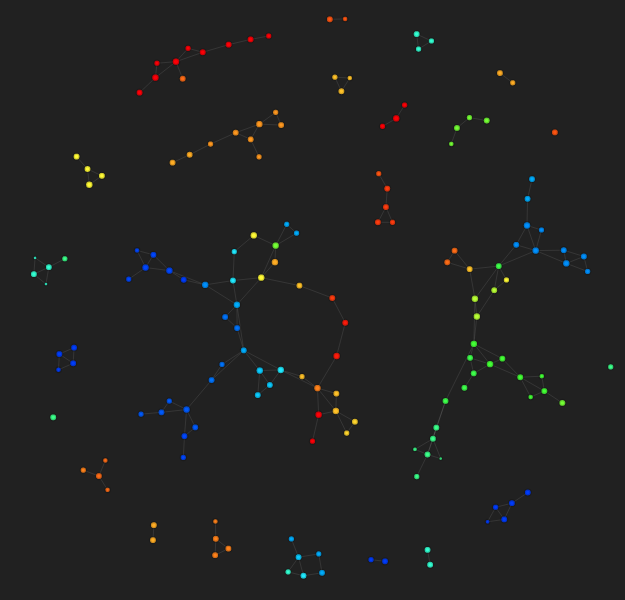}}
\caption{Mapper graphs obtained using Keplermapper \cite{hendrik_jacob_van_veen_2017_1002378} via $2$D $t$-SNE lens on UCI \cite{UCI} Dota 2 Games Results Data Set restricted to games involving playable character \textit{Dazzle}.}
\label{fig:multimapper motivation}
\end{figure*}

It is possible, and we have seen in experiments, that the same cover might not be ideal for every part of the data. For example, in Fig. \ref{fig:multimapper motivation}, the coarsest refinement hides local structure that appears at finer scales, but the graph also begins to break into more components, which obscures global relationships between the data in each component.\\
To solve this issue, we have developed a technique which takes the locally best scale for each region of the data, and glues them together to represent the entire data more accurately than a single global scale.\\
For scale selection, we use the notion of Multiscale Mapper. However, the varying scales in Multiscale Mapper is applied globally, which means that two regions with different density would not be appropriately represented at any one level. Therefore we propose a Mapper graph which is an amalgamation of Mapper graphs computed on various regions of interest. We apply Multiscale Mapper on each region of interest as a way of selecting scales.
It is a different question to identify these regions and the appropriate scale of its cover in the first place -- we discuss a possible approach in \ref{sec:det bad} using the idea of Multiscale Mapper.\\
To consistently combine the locally suitable Mappers, we implement our idea as a repeated rescaling of selected regions. Intuitively it can be understood as magnification/compression of a Mapper in relatively sparser/denser regions respectively. Our first approach slices relevant bins of the cover to obtain smaller bins; our second approach is more sophisticated and uses a nerve-like computation to glue together various locally suitable Mappers. The latter has the advantage of being compatible with all covering schemes, including the brick-like cover proposed in \ref{sec:brickcover}.

\subsection{Local Refinement via Cover Slicing}\label{multi app 1}
\begin{figure}[b]
\centering
\includegraphics[width=0.4\linewidth]{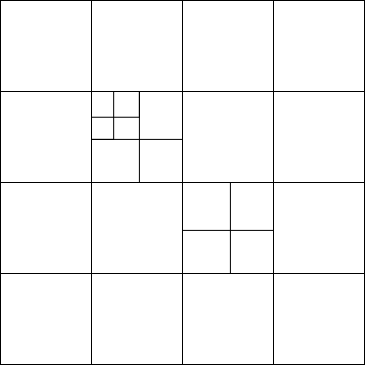}
\caption{Original bins sliced to obtain smaller bins. As shown here, the process can be repeated.}
\label{fig:sliced cover}
\end{figure}

The covering process can be broken into two parts, deciding the type of partition and then overlap percentage. Our aim in this section is to modify the underlying partition in a manner that the required regions are covered by smaller bins than before.
Once we have identified regions in the Mapper $M$ that are to be magnified, we perform the steps below. It is specifically illustrated using cuboidal bins, but can be generalized to other shapes by redefining the \textit{chopping} in Step 3.
\begin{enumerate}
\item Let S be the set of nodes we wish to magnify; $\tilde{X}$ be the region of the original data corresponding to these nodes, i.e. 
\begin{equation}
    \tilde{X} = \cup_{w \in S}C_w \label{eq7}
\end{equation}
where $C_w$ is the cluster corresponding to a node $w$.
\item Now we look at the image of $\tilde{X}$ under the original lens function $f$, i.e. 
\begin{equation}
\tilde{Z} = f(\tilde{X})
\end{equation}
\item $\Part$ be the partition of the image, $Z = f(X)$, that gave rise to the Mapper $M$. Then define a subset:
\begin{equation}
    \Part' = \{P \in \Part : P \cap \tilde{Z} \neq \emptyset\}
\end{equation}
i.e. those parts of the partition which contain some part of $\tilde{Z}$.
\item Obtain a refinement $\Part''$ of $\Part$ by slicing each box along each dimension. For example, a 1D interval would be halved; a 2D box would be sliced into four equal 2D boxes, as shown in Fig \ref{fig:sliced cover}. This is a \textit{refinement by $2$} -- slicing into $m$ pieces along each axis would be a \textit{refinement by $m$}.
\item Obtain a new Mapper $\tilde{M}$ with the same lens function, but with a cover built from: \begin{equation}\tilde{\Part} = (\Part \setminus \Part')\cup \Part''\end{equation}
\end{enumerate}
$\tilde{M}$ thus shows the data corresponding to $\tilde{Z}$ $M$, and the remaining parts of the data at the same scale as $M$. This process can be repeated at various regions of the data to view each region at a suitable scale.
\subsection{Handling Degeneracy via Multimapper}\label{multi app 2}
A major drawback of the previous approach is degeneracy of $f$ i.e. $f$ may map distant parts of $X$ very closely, in the image $Z$. Hence, if we magnify a region $\tilde{X}$ in the above method, since we go via its image $\tilde{Z}$, we would actually magnify $f^{-1}(\tilde{Z}) = f^{-1}(f(\tilde{X}))$, which is potentially a larger region and having other parts distant from $\tilde{X}$. Hence undesirable parts of the Mapper would be magnified as well. We wish to avoid this effect of $f$; but we want to retain the convenience of constructing the cover via pullback under $f$. Moreover, we would want to not only \textit{zoom in} on denser parts, but \textit{zoom out} on sparser parts that might have shattered. \\
Our one-shot solution to these requirements is the notion of Multimapper. Multimapper is constructed as follows:
\begin{figure*}
\centering
\subfloat[Bin size = $1/5$\textsuperscript{th} of image diameter. In white : nodes selected for magnification]{\includegraphics[width=0.27\linewidth,height=4.5cm]{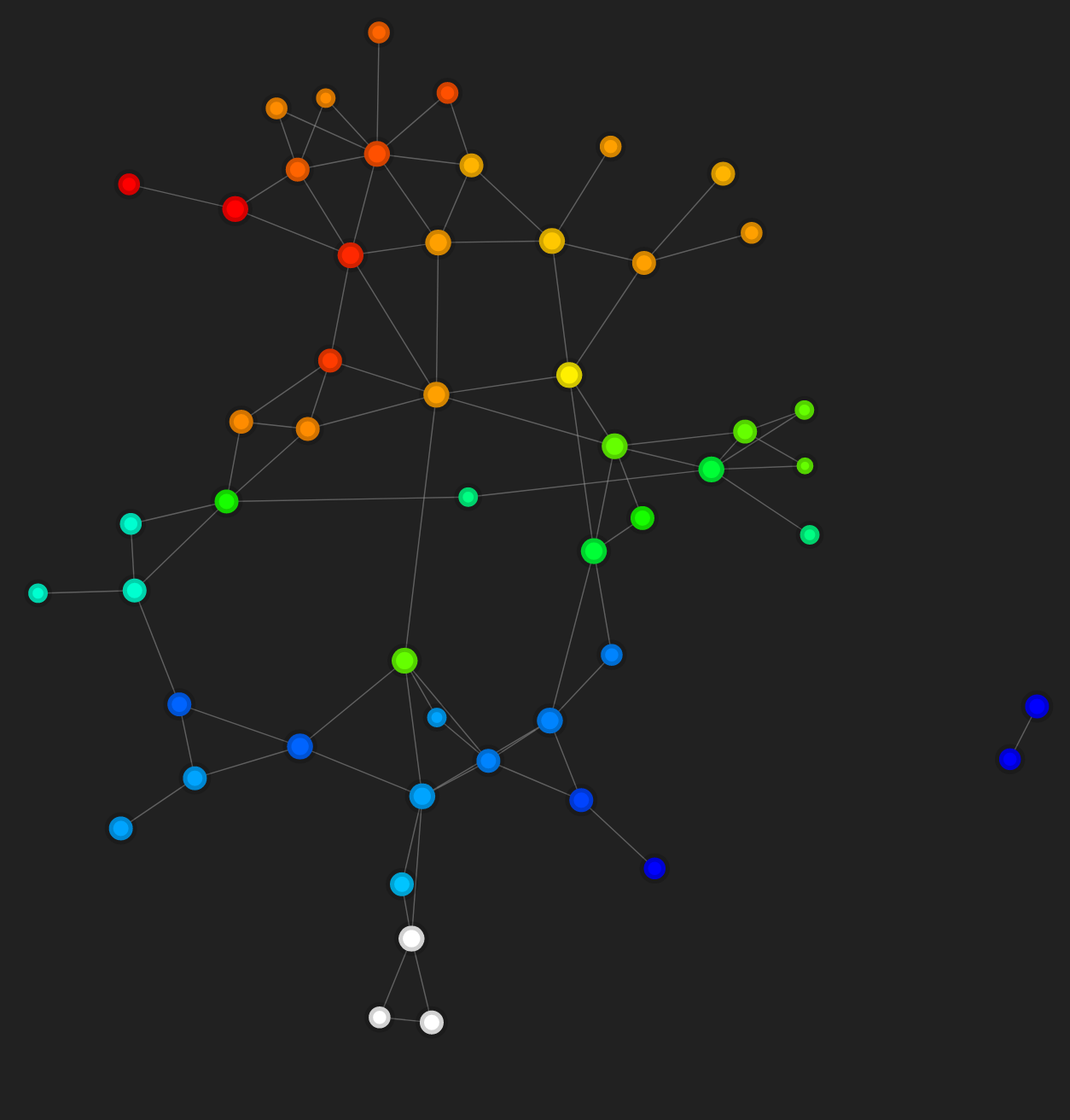}} \qquad
\subfloat[magnification of selected region. Bin size in selected region = $1/10$\textsuperscript{th} of image diameter.]{\includegraphics[width=0.27\linewidth,height=4.5cm]{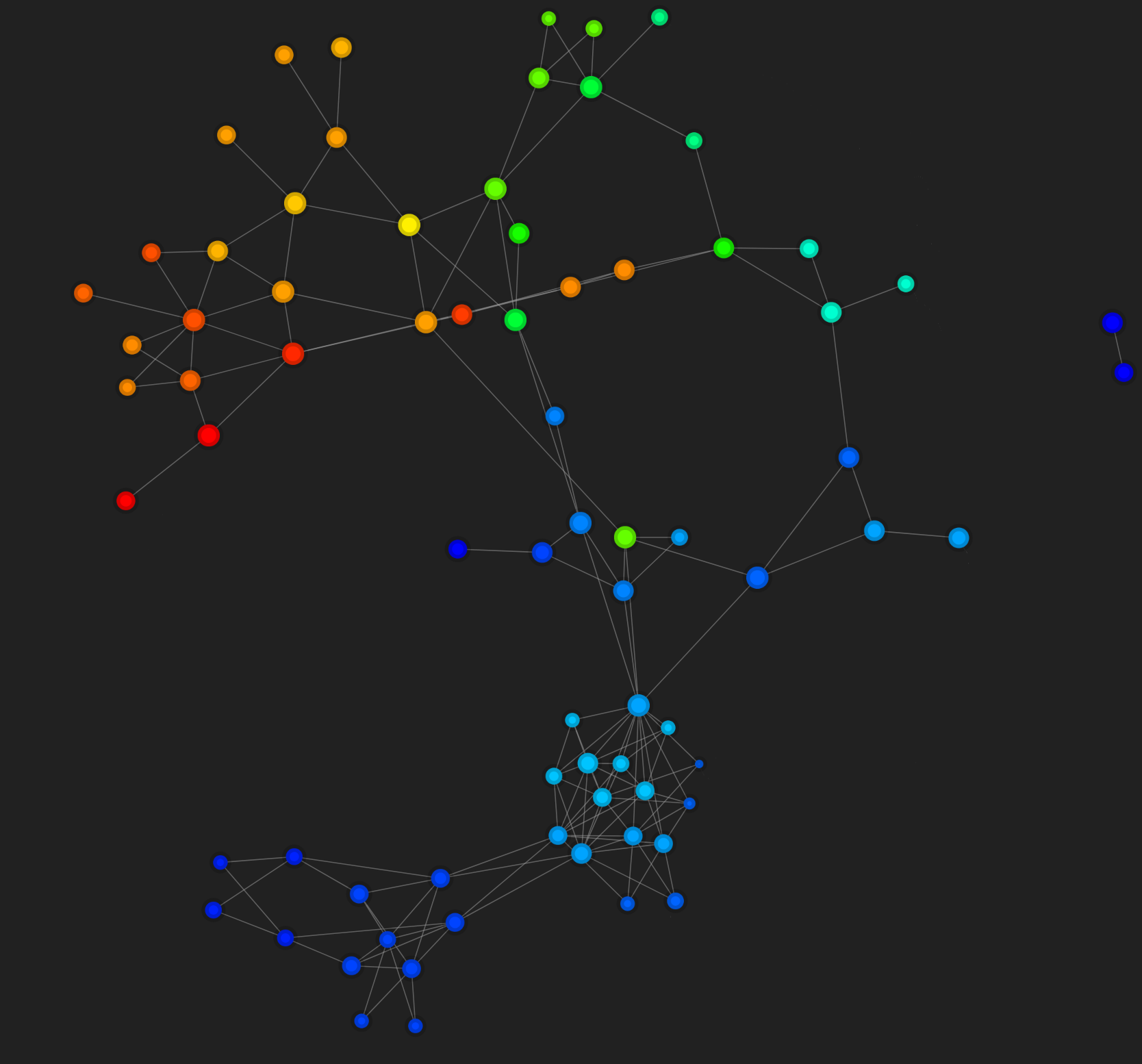}} \qquad
\subfloat[Final Multimapper graph. Bin size for magnification = $1/10$\textsuperscript{th} of image diameter. ]{\includegraphics[width=0.27\linewidth,height=4.5cm]{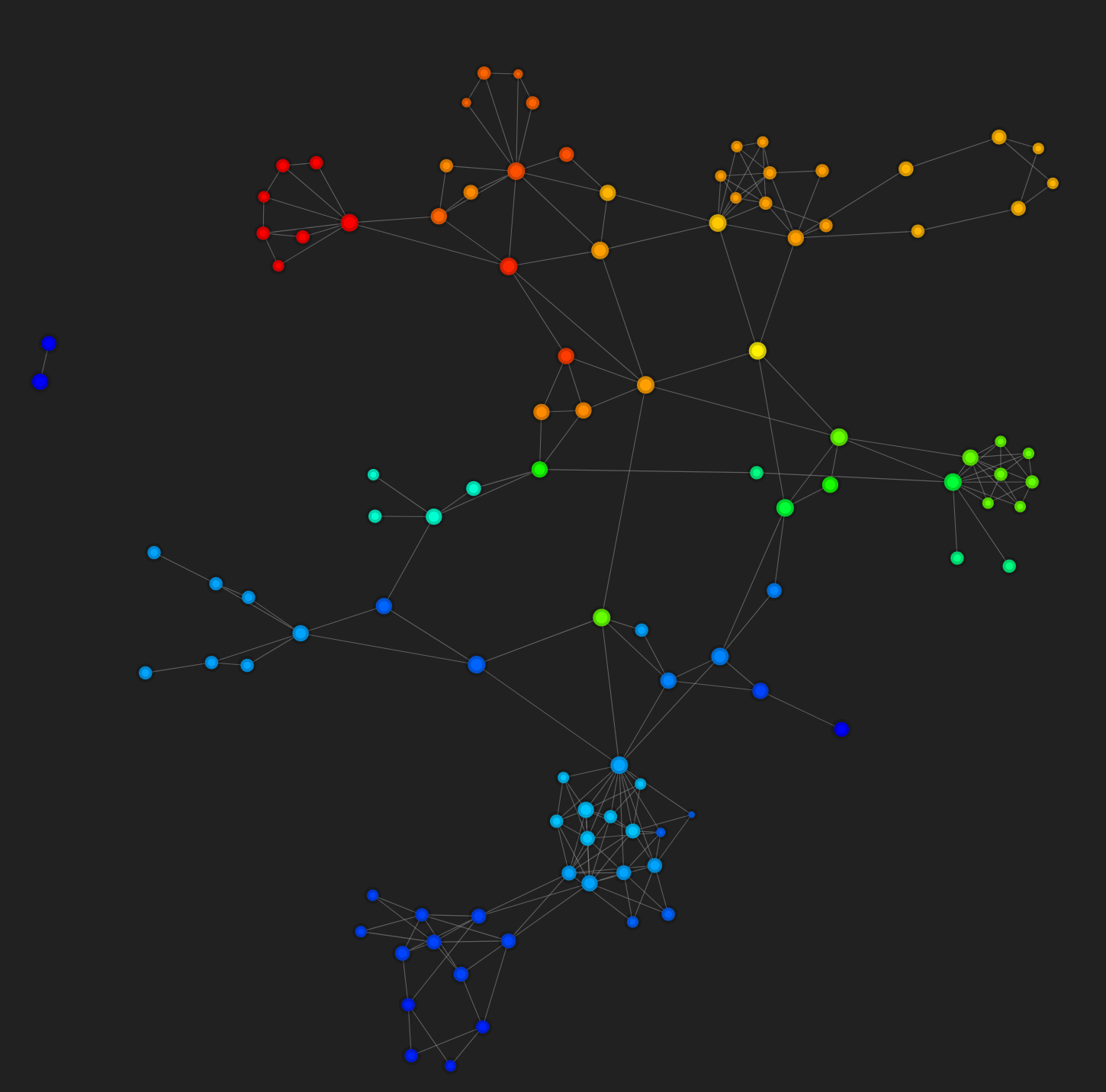}}
\caption{Mapper graphs of the Dota 2 Data Set from Fig. \ref{fig:multimapper motivation}, using Keplermapper \cite{hendrik_jacob_van_veen_2017_1002378} with $2$D $t$-SNE filter (left); Multimapper graphs (b) and (c) which reveals new geometry.}
\label{fig:magnification}
\end{figure*}
\begin{enumerate}
\item As before, given a Mapper graph $M$, with nodes $V$ and some nodes $S$ to be magnified, we identify the corresponding data subset and image subset:
\begin{equation}
\begin{aligned}
\tilde{X} &= \cup_{w \in S}C_w \\
\tilde{Z} &= f(\tilde{X})
\end{aligned}
\end{equation}
where $C_w$ is the cluster corresponding to a node $w$.
\item Let $\Clus$ be the set of clusters corresponding to the nodes of $M$. The region we want to `preserve at original scale', i.e. \emph{not} magnify, is: \begin{equation}X' = X \setminus \tilde{X}\end{equation} From this we discard those clusters which lie entirely in $\tilde{X}$ to obtain: \begin{equation}
\begin{aligned}
&\Clus' = \{C \in \Clus : C \cap X' \neq \emptyset\} \\
&such\;that:\;\cup_{C \in \Clus'}C \supseteq X'
\end{aligned}
\end{equation} 
\item We define a new cover on $\tilde{Z}$ as per our choice and requirement. Unlike in \ref{multi app 1}, this need not be a restriction of the old cover of $Z$.
\item We cluster within the inverse images of each new bin, to obtain a new set of clusters $\tilde{\Clus}$, such that: \begin{equation}\cup_{C \in \tilde{\Clus}}C = \tilde{X}\end{equation} Effectively, it is again a Mapper construction restricted on $\tilde{X}$ with a new cover on $\tilde{Z}$.
\item Now we have obtained an overall set of clusters 
\begin{equation}
    \begin{aligned}
    &\hat{\Clus} = \Clus' \cup \tilde{\Clus}\\
    &such\;that:\;\cup_{C \in \hat{\Clus}}C = X
    \end{aligned}
\end{equation} We compute a new nerve $\hat{M}$ according to $\hat{\Clus}$: the nodes correspond to clusters in $\hat{\Clus}$, and a $k$-simplex is added for every $k+1$ clusters of $\hat{\Clus}$ that have a simultaneous intersection.\\
This nerve-like computation ensures that $\hat{M}$:
\begin{itemize}
\item matches $M$ on $V \setminus S$, via nerve computation on $\Clus'$
\item replaces the induced sub-complex on $S$ with a copy of $\tilde{M}$, via nerve computation on $\tilde{\Clus}$
\item \textit{glues} these two parts in a manner faithful to the topology of $X$, by imitating the nerve construction on simultaneous intersections involving both $\Clus'$ and $\tilde{\Clus}$
\end{itemize}
\end{enumerate}
The above process can be repeated at various regions, with locally suitable choices of covering scheme, bin size, and clustering algorithm. The resulting structure is a Mapper-like simplicial complex which we call \emph{Multimapper}.\\
Conceptually, Multimapper breaks up the original point cloud into subsets that may or may not intersect. For each region, it computes the Mapper that is of a suitable scale for that region.\\
\begin{defn}[Multimapper]
Given a finite collection of subsets $\mathcal{X} \subset 2^X$ such that $\cup_{Y \in \mathcal{X}}Y=X$ and covers $\{\Cover_Y\}_{Y \in \mathcal{X}}$ such that $\forall Y \in \mathcal{X}, \Cover_Y \supseteq Y$, we can define the corresponding \emph{Multimapper}: \begin{equation}\Multi(\mathcal{X}, \{\Cover_Y\}_{Y \in \mathcal{X}}) = N(\cup_{Y \in \mathcal{X}}f^*(\Cover_Y))\end{equation}
\end{defn}
This piecewise approach makes Multimapper quite flexible, we can freely choose different local Mappers for different regions, while retaining the global relationships between them. Fig. \ref{fig:magnification}(b) shows an example of applying Multimapper on a specific region of dataset. Fig. \ref{fig:magnification}(c) shows the graph obtained by applying Multimapper repeatedly on separate regions and creating a single visualization by gluing them all. We can see it reveals new structure in the magnified regions as well as prevents the shattering of graph which was happening in case of Mapper with similar bin size as can be seen in Fig. \ref{fig:multimapper motivation}(c). \\
Via the stability results in \cite{dey2016multiscale}, we further know that Multimapper is locally as stable as Multiscale Mapper; and of course, since it after all arises from a cover, it is globally as stable as Mapper.
\section{Detecting a Bad Mapper}\label{sec:det bad}
\cite{singh2007topological} refer to the \emph{Reeb graph} as a geometric representation suitable for obtaining information directly, and introduce Mapper as a generalization of Reeb graph. Later, \cite{carri19statistical} have shown that in the ideal setup, Mapper with a $1$-dimensional lens function statistically converges to the Reeb graph of the original space under the lens function. A general convergence result of Mapper to a generalization of Reeb graph called the \emph{Reeb space} has been conjectured, and \cite{CatReeb} have studied a category-theoretic version of this relationship.\\
The ideal convergence would not occur in case of real data since real data is a discrete point cloud. However, we can measure the correctness of a Mapper via its closeness to the Reeb space. We can thus reasonably demand that a good enough Mapper of a discrete point cloud, under
a particular lens function $f$, should be similar in shape to the Reeb space (under $f$) of the connected paracompact space it approximates.
\begin{defn}[Reeb space]\label{reeb space}
Given a continuous map $f:X \rightarrow Y$ between topological spaces $X$ and $Y$ , the \emph{Reeb space} $R_f(X)$ of $X$ with respect to $f$ is $X/\sim$ where the equivalence relation $\sim$ on $X$ is defined as $p \sim q$ iff $p$ and $q$ lie in the same connected
component of $f^{-1}(c)$ for some $c \in Y$.
\end{defn}
When $Y = \R^n$, the connected components are same as path connected components \cite{sutherland2009introduction}, and this is sufficient for our real-world setting. When $Y = \R$, the Reeb space is called the \emph{Reeb graph}.

The Mapper algorithm is constructed using the Nerve Theorem, and we have characterized its goodness by its closeness to the Reeb space. Thus, to partially characterize bad Mappers i.e. Mappers far from the Reeb space, we can try to find regions of the Mapper where the contractibility hypotheses of Nerve Theorem is violated.
\\Checking for contractibility, however, is complicated in a real world setting, especially since our actual space is a point cloud.\\
To adapt our method to the real world, let us recall the association between the continuous and discrete ideas:
\begin{itemize}
\item Continuous $\leftrightarrow$ Discrete
\item Paracompact space $\leftrightarrow$ Point cloud
\item Path connected components $\leftrightarrow$ Clusters
\item Nerve $\leftrightarrow$ Mapper
\end{itemize}

We know that if a space has more than one path connected component, it cannot be contractible -- shrinking two separate pieces to the same point would require shrinking across the gap between them, which violates continuity. Translated to the point cloud setting, this means that if a data subset has multiple clusters, the corresponding space cannot be contractible. Hence, a sufficient condition for non-contractibility can be checked using the following general characterization: \\
Given a Mapper on a set of nodes $V$
\begin{equation}\exists \sigma \in M, \sigma = \sigma(S), S \subseteq V, \beta_0({\cap_{v \in S}C_v}) > 1\end{equation}
where:
\begin{itemize}
\item $\sigma(S)$ is the simplex on the vertices $S$
\item For some topological space $A$, $\beta_0(A)$ is the number of connected components in $A$
\item $C_v$ is the cluster corresponding to the node $v$ of the nerve.
\end{itemize}
This motivates our approach, which we illustrate via $1$-simplices (edges). However, the same technique can be iterated over all simplices in $M$. The most naive approach of identifying components in discrete setting is via any known clustering algorithms. We propose such a method next, followed by a modification on it that is independent of clustering algorithms.
\subsection*{Clustering-Dependent Version}\label{nerve app 1}
\textbf{Procedure:} Given a Mapper $M$, for each edge $(u, v)$, $C_u \cap C_v \neq \emptyset$, we cluster within $C_u \cap C_v$ using a clustering algorithm like DBSCAN which does not fix the number of clusters \textit{a priori}. If more than one cluster is obtained, we report it as a violation. Finding even a single violation is sufficient for the Mapper to be classified as \textit{bad}.\\
\textbf{Explanation:} If an edge $e = (u, v)$ leads to a violation, it means that the continuous space $C_u \cap C_v$ approximates has more than one path connected component. Hence the Reeb space of $C_u \cap C_v$ , obtained by collapsing each path connected component, will also have more than one connected component. But the Mapper restricted to $C_u \cap C_v$ is precisely the edge $e$, which is a single connected component. Hence in the region of data corresponding to $C_u \cap C_v$, the Mapper deviates in shape from the Reeb space.

\subsection{Our Algorithm: Clustering-Independent Version}\label{nerve app 2}
In the above naive approach, we depend on clustering to approximate path connected components. Thus we are constrained by the choice of clustering algorithm and its parameters, and must optimize this on a case-by-case basis as these are not generalizable to any dataset. To remove this dependency, we propose a method of approximating path connected components that is independent of clustering.\\
The well-known TDA method of \emph{persistence}, which is usually applied directly on the point cloud \cite{ghrist2008barcodes}, can be translated to the Mapper setting via Multiscale Mapper -- \cite{dey2016multiscale} have given an algorithm that, given a tower of covers, computes the \emph{persistence diagram} of the resulting Multiscale Mapper. They have also provided an approximate computation suitable for the discrete point cloud setting.\\
Hence, given $\Tower = \{\Cover_\eps\}_\eps$, a Multiscale Mapper $\MM(\Tower)$ and a void $H$ that appears in some $M(\Cover_\eps) \in \MM(\Tower)$
\begin{enumerate}
\item $birth(H) = \textrm{min}\{\eps : H \text{ appears in }M(\Cover_\eps)\}$
\item $death(H) = \textrm{max}\{\eps : H \text{ appears in }M(\Cover_\eps)\}$
\end{enumerate}
Hence, for every $0$-void, i.e. connected component, $1$-void, i.e. circular hole and so on, that appears in $\MM(\Tower)$, we get a birth-death pair, a range of scales at which it is visible in the corresponding Mappers.\\
We consider a topological feature, e.g. a void, to be truly present in the shape of the data if it remains, or is \emph{persistent}, for a large range of scales. If $\beta_m^\eps$ is the number of $m$-voids at scale $\eps$, we remove the noisy features and calculate the \textit{true} number of $m$-voids $\beta_m$ i.e. the persistent ones as: 
\begin{equation}
\begin{aligned}
\phi_m : \N & \rightarrow \N \\
\beta & \mapsto \abs{\{\eps : \beta_m^\eps = \beta\}}\\
\beta_m &= \textrm{arg max} (\phi_m)
\end{aligned}
\end{equation}
Since $\beta_0$ is the number of connected components, persistence on Mapper via Multiscale Mapper gives us a cluster-independent way to compute the number of connected components in the data. This gives us the following procedure for detecting violations:\\
Given a Mapper $M(X,\Cover,f)$, for each edge $(u,v)$, we construct a tower of covers on $C_u \cap C_v$, beginning from $f^*(\Cover)\restrict{C_u\cap C_v}$ and decreasing the scale of cover up to a threshold of refining $f^*(\Cover)\restrict{C_u\cap C_v}$ by $2$.\\
Thus in practice, where $\Cover\restrict{f(C_u\cap C_v)}$ has cuboidal\footnote{Cuboidal bins may be arranged in the standard way or as suggested in \ref{sec:brickcover}.} bins of diameter $\eps_0$ and $\Tower'$ must be finite, we define covers of the form $\Cover_\eps$ of $f(X)$ to have the same partition rule as $\Cover$, but with bin size $\eps$. Thus the tower of covers of $f(C_u\cap C_v)$ is:
\begin{equation}\Tower'=\{\Cover_\eps\}_{\eps \in \N, 0.5\eps_0 \leq \eps \leq \eps_0}\end{equation}
From this we obtain the pullback tower $f^*(\Tower')$ of covers of $C_u\cap C_v$. On this, using persistence via Multiscale Mapper as illustrated in \cite{dey2016multiscale}, we can compute $\beta_0$, the persistence and hence \textit{true} number of components. If $\beta_0 > 1$, we report a violation. To increase efficiency in implementation, we will compute only zeroth dimension persistence diagram.


\section{Reducing Obscured Information via Brick-like Cover}\label{sec:brickcover}
Construction of a cover from the implementation perspective for the Mapper algorithm can be conceptualized in two steps:

\begin{enumerate}
\item Construction of a partition

\item Growing each piece of the partition to introduce overlap, hence obtaining bins.
\end{enumerate}
The most well-known, standard box-like cover partitions the image space into $n$-cuboids, where $n$ is the projected dimension. Hence with a $1$-dimensional lens, we get line segments; with a $2$-dimensional lens, we get rectangular boxes; and so on. In such a setup, $2^n$ bins can intersect where $2^n$ pieces of the partition meet; hence simplices in the nerve can have dimension up to $2^n$. However, simplices of dimension $>2$ are difficult to visualize. Because of this, the visualization we create will be truncated at $2$- or $3$-simplices, and not represent the complete information contained in the Mapper.\\
Hence, it would be desirable to build covers such that the visualization represents all the topological information contained in Mapper through simplices of visualizable dimensions. This requires us to explore non-cuboidal lattices for the partition. As suggested in \cite{singh2007topological}, hexagonal lattices are a natural choice -- in $2$D, for example, hexagonal bins can intersect only $3$ at a time; hence the Mapper would have at most $2$-simplices.

\subsection*{The Brick Cover}

Constructing bins over a hexagonal lattice is computationally difficult even in $2$ dimensions. Our proposed $2$D covering is more efficient and achieves the same goal with a few realistic constraints. Our partition of the $2$D plane uses offset rectangles, as in a brick wall.
The underlying vertices are still a hexagonal lattice, which gives us an advantage over the usual rectangles: at most $3$ bricks can meet at a vertex. At overlap below $50\%$, this translates to a maximum of $3$-fold intersection of bins. Given the sparse nature of high-dimensional data, overlap below $50\%$ is a reasonable notion of nearness.\\
The idea of using non-cuboidal bins can be extended to higher dimensions: for example, cuboidal lattice in $3$D gives up to $8$-fold intersections, while hexagonal lattice in $3$D gives only up to $6$-fold intersections.
\begin{figure}[t]
\centering
  \includegraphics[width=0.5\linewidth]{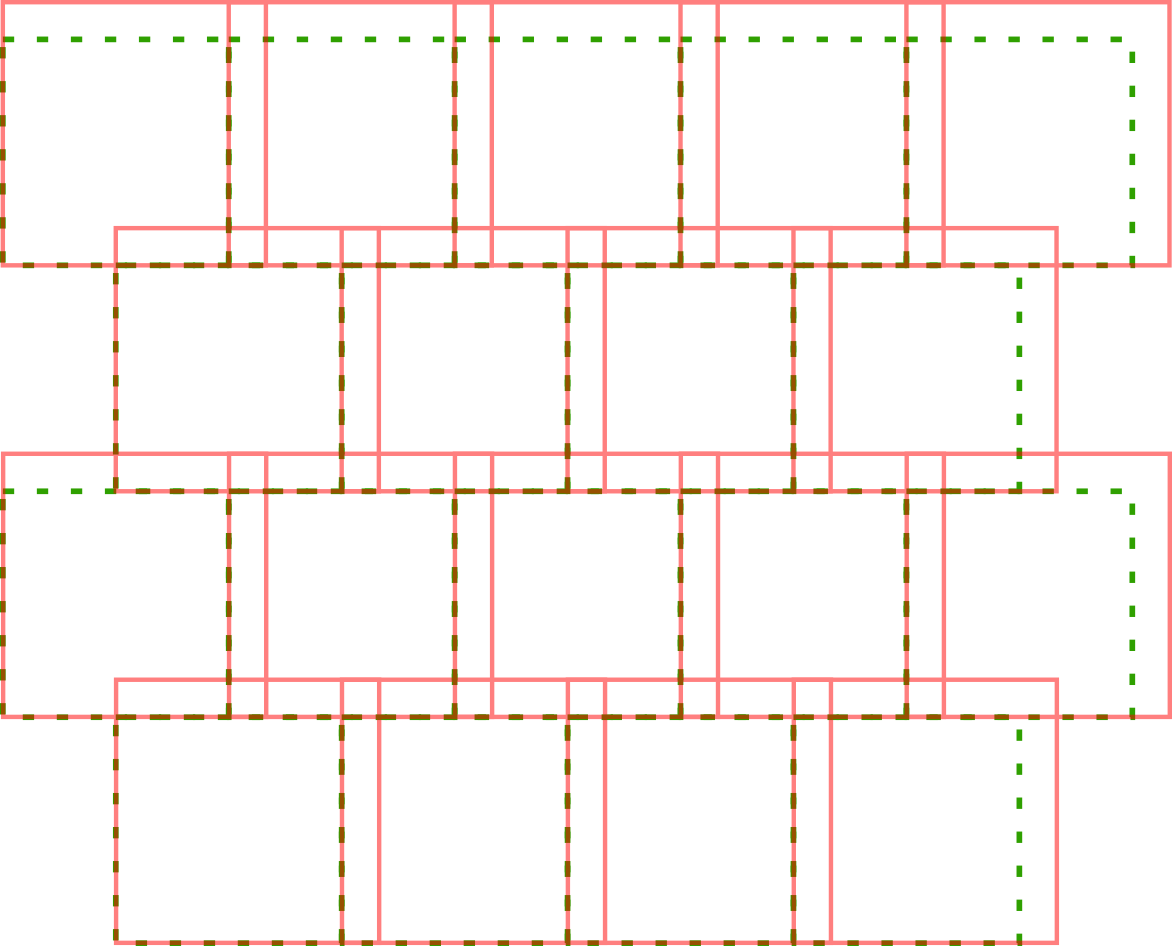}
  \caption{Bins built on brick-like mesh. Notice how at most $3$ bins intersect at once.}
  \label{fig:brickcover}
\end{figure}
We \emph{claim that our proposed brick-like cover gives the lowest possible maximum simultaneous intersections among all $2$D covers that use cuboidal bins}. This is because:
\begin{enumerate}
\item For the underlying partition of a cuboidal cover, let a `mesh point' be a point at which some pieces of the partition meet, and at least one piece has that point as a vertex. At a mesh point:
\begin{itemize}
\item The total angle contributed by pieces incident on it must equal $2\pi$.
\item A piece having the mesh point as its own vertex, contributes $\pi/2$, other surrounding pieces contribute $\pi$. Hence $2$ pieces cannot form a mesh point, since neither piece could have the said point as its own vertex.
\item Thus, at a minimum, we are forced to build $2\pi$ as $\pi+\pi/2+\pi/2$, as in a brick-like mesh.
\end{itemize}
\item We are left to choose the `offset', i.e. by how much the successive rows of pieces are shifted from each other. We assume that the overlap is added to the top right of the underlying pieces. Then, between two successive rows:
\begin{itemize}
\item If the top row is shifted to the right by $p\%$, then an overlap greater than $p\%$ would lead to $4$-intersections. So we must maximize $p$. But this is the same as the bottom row being shifted to the left by $(100-p)\%$, so an overlap greater than $(100-p)\%$ would also cause $4$-intersections. So we must maximize $(100-p)$.
\item Symmetrically, when the top row is shifted to the left by $p$, we again need a simultaneous maximization of $p$ and $(100-p)$.
\item Hence we must choose $p = (100-p) = 50$, which gives us the proposed cover.
\end{itemize}
\end{enumerate}

\section{Conclusion}\label{sec:conclusion}
In this paper, we proposed improvements upon the existing Mapper algorithm solving many of its shortcomings partially. We have given a partial characterization of undesirable outputs of the Mapper algorithm and proposed a flexible method that corrects the choice of scale locally in a manner sensitive to the density of various data subsets. In all these methods, we have retained the \emph{unsupervised} nature and stability of Mapper. Moreover, replacing the standard covering scheme with our brick-like cover reveals more topological information in a visualizable way.
Our methods produce a visualization that is more true to the actual shape of data, via the Reeb space characterization, than the standard Mapper. Our contributions pave the path towards an automatic one-shot Mapper output which is the best visualization in terms of being close to the topological structure of data without any need of manual parameter optimization. This improves the efficacy of its applications in analysis and visualization of high dimensional big data.\\ 
An interesting direction of future work that we mean to pursue is to study the relationship between successive applications of the Multimapper algorithm.
Moreover, our method to detect deviations from Reeb space in Mapper via violations of Nerve Theorem hypothesis can be more powerfully implemented if a discrete analog of \textit{contractibility} was reasonably defined, similar to how clusters are used to approximate connected components. Finally, a characterization of the best cover possible in any given dimension for the given data would be useful.

\bibliographystyle{apalike}
\bibliography{references}

\begin{thebibliography}{}

\bibitem[Ayasdi, 2015]{AyasdiHealth}
Ayasdi, I. (2015).
\newblock Clinical variation management with advanced analytics.

\bibitem[Ayasdi, 2016]{AyasdiFinance}
Ayasdi, I. (2016).
\newblock Machine intelligence for financial services.

\bibitem[Ayasdi, 2017]{AyasdiML}
Ayasdi, I. (2017).
\newblock Tda and machine learning: Better together.

\bibitem[Ayasdi, 2018]{AyasdiMain}
Ayasdi, I. (2018).
\newblock Understanding ayasdi: What we do, how we do it, why we do it.

\bibitem[Beckham, 2012]{AyasdiNBA}
Beckham, J. (2012).
\newblock Analytics reveal 13 new basketball positions.

\bibitem[Carriere et~al., 2018]{carri19statistical}
Carriere, M., Michel, B., and Oudot, S. (2018).
\newblock Statistical analysis and parameter selection for mapper.
\newblock {\em The Journal of Machine Learning Research}, 19(1):478--516.

\bibitem[Dey et~al., 2016]{dey2016multiscale}
Dey, T.~K., M{\'e}moli, F., and Wang, Y. (2016).
\newblock Multiscale mapper: Topological summarization via codomain covers.
\newblock In {\em Proceedings of the twenty-seventh annual acm-siam symposium
  on discrete algorithms}, pages 997--1013. SIAM.

\bibitem[Dheeru and Karra~Taniskidou, 2017]{UCI}
Dheeru, D. and Karra~Taniskidou, E. (2017).
\newblock {UCI} machine learning repository.

\bibitem[Ghrist, 2008]{ghrist2008barcodes}
Ghrist, R. (2008).
\newblock Barcodes: the persistent topology of data.
\newblock {\em Bulletin of the American Mathematical Society}, 45(1):61--75.

\bibitem[Hatcher, 2002]{Hatcher}
Hatcher, A. (2002).
\newblock {\em Algebraic topology}.
\newblock Cambridge University Press, Cambridge.

\bibitem[Kamruzzaman et~al., 2016]{kamruzzaman2016characterizing}
Kamruzzaman, M., Kalyanaraman, A., and Krishnamoorthy, B. (2016).
\newblock Characterizing the role of environment on phenotypic traits using
  topological data analysis.
\newblock In {\em Proceedings of the 7th ACM International Conference on
  Bioinformatics, Computational Biology, and Health Informatics}, pages
  487--488. ACM.

\bibitem[Munch and Wang, 2016]{CatReeb}
Munch, E. and Wang, B. (2016).
\newblock Convergence between categorical representations of reeb space and
  mapper.
\newblock In {\em 32nd International Symposium on Computational Geometry, SoCG
  2016, June 14-18, 2016, Boston, MA, {USA}}, pages 53:1--53:16.

\bibitem[Pearson, 1901]{pearson1901principal}
Pearson, K. (1901).
\newblock Principal components analysis.
\newblock {\em The London, Edinburgh, and Dublin Philosophical Magazine and
  Journal of Science}, 6(2):559.

\bibitem[Singh et~al., 2007]{singh2007topological}
Singh, G., M{\'e}moli, F., and Carlsson, G.~E. (2007).
\newblock Topological methods for the analysis of high dimensional data sets
  and 3d object recognition.
\newblock In {\em SPBG}, pages 91--100.

\bibitem[Steen et~al., 1978]{steen1978counterexamples}
Steen, L.~A., Seebach, J.~A., and Steen, L.~A. (1978).
\newblock {\em Counterexamples in topology}, volume~18.
\newblock Springer.

\bibitem[Sutherland, 2009]{sutherland2009introduction}
Sutherland, W.~A. (2009).
\newblock {\em Introduction to metric and topological spaces}.
\newblock Oxford University Press.

\bibitem[van Veen and Saul, 2017]{hendrik_jacob_van_veen_2017_1002378}
van Veen, H.~J. and Saul, N. (2017).
\newblock Keplermapper.

\bibitem[Vejdemo-Johansson et~al., 2012]{vejdemo2012topology}
Vejdemo-Johansson, M., Carlsson, G., Lum, P.~Y., Lehman, A., Singh, G., and
  Ishkhanov, T. (2012).
\newblock The topology of politics: voting connectivity in the us house of
  representatives.
\newblock In {\em NIPS 2012 Workshop on Algebraic Topology and Machine
  Learning}.

\end{thebibliography}

\end{document}